\documentclass[12pt]{amsart}

\usepackage{hyperref}
\usepackage{fullpage}
\usepackage{amsmath}
\usepackage{amsthm}
\usepackage{mathrsfs}
\usepackage{amssymb}
\usepackage{comment}
\usepackage{fancyhdr}
\newtheorem{lemma}{Lemma}
\newtheorem{theorem}{Theorem}

\newtheorem{prop}{Proposition}

\newtheorem{remark}{Remark}

\newtheorem{corollary}{Corollary}
\newtheorem{definition}{Definition}
\newtheorem{observation}{Observation}
\newcommand{\Z}{\mathbb{Z}}
\newcommand{\C}{\mathbb{C}}

\begin{document}

\title{Proof of a conjecture of Guy on class numbers}
\author{Lynn Chua}
\address{Department of Mathematics, Massachusetts Institute of Technology, 77 Massachusetts Avenue
Cambridge, MA, 02139}
\email{chualynn@mit.edu}
\author{Benjamin Gunby}
\address{Department of Mathematics, Massachusetts Institute of Technology, 77 Massachusetts Avenue
Cambridge, MA, 02139}
\email{bgunby@mit.edu}
\author{Soohyun Park}
\address{Department of Mathematics, Massachusetts Institute of Technology, 3 Ames St., Cambridge, MA, 02139}
\email{soopark@mit.edu}
\author{Allen Yuan}
\address{Department of Mathematics, Harvard University, 1 Oxford St., Cambridge, MA, 02138}
\email{allenyua@gmail.com}

\maketitle

\begin{abstract}
It is well known that for any prime $p\equiv 3$ (mod $4$), the class numbers of the quadratic fields $\mathbb{Q}(\sqrt{p})$ and $\mathbb{Q}(\sqrt{-p})$, $h(p)$ and $h(-p)$ respectively, are odd.  It is natural to ask whether there is a formula for $h(p)/h(-p)$ modulo powers of $2$.  We show the formula $h(p) \equiv h(-p) m(p)$ (mod $16$), where $m(p)$ is an integer defined using the ``negative" continued fraction expansion of $\sqrt{p}$. Our result solves a conjecture of Richard Guy.
\end{abstract}

\section{Introduction} \label{sectintro}
In this paper, $p$ denotes a prime which is congruent to 3 modulo 4, and $h(p)$ (respectively $h(-p)$) denotes the class number of the quadratic field $\mathbb{Q}(\sqrt{p})$ (respectively $\mathbb{Q}(\sqrt{-p})$).  Class numbers are a deep subject which has been extensively studied.  For example, the Scholz reflection principle \cite{Pierce} gives information about the 3-part of class numbers of corresponding real and imaginary quadratic fields. In the spirit of such results, we will consider the relationship between $h(p)$ and $h(-p)$ modulo powers of two. Some results in this area are already known; for example, it is known that $h(p)$ and $h(-p)$ are odd (see \cite{Brown1}, \cite{Brown2}). In addition, we have the following result modulo $4$.

\begin{theorem}[Mordell, \cite{Mordell}]
Let $p > 3$. Then, we have
\begin{equation*}
h(-p) \equiv
\begin{cases}
1 \, \pmod{4}\, & \text{if } \left(\frac{p - 1}{2}\right)! \equiv -1 \pmod{p}\,, \\
3 \, \pmod{4}\, & \text{if } \left(\frac{p - 1}{2}\right)! \equiv 1 \pmod{p}\,.
\end{cases} \,
\end{equation*}
\end{theorem}

Williams (\cite{Will79}, \cite{Will82}) obtains results on the relationship between $h(p)$ and $h(-p)$ modulo 4 and 8. 

\begin{theorem}[Williams, \cite{Will79}, \cite{Will82}] \label{Willthm}
Let $p > 3$ and let $T + U\sqrt{p}$ be the fundamental unit of $\mathbb{Q}(\sqrt{p})$. 
Then 
\begin{eqnarray*} h(-p) &\equiv& h(p)(U + 2) \pmod 4\,, \\
 h(-p) &\equiv& h(p)\left(2 + pU - 2\left(\frac{T}{U}\right)\right) \pmod 8\,.\end{eqnarray*}
\end{theorem}

Motivated by a conjecture of Richard Guy \cite{Guy}, we study congruences between $h(p)$ and $h(-p)$ modulo $16$. Let the following be the ``negative" continued fraction expansion of $\sqrt{p}$
\begin{equation} \label{contfracp}
\sqrt{p} = b_0 - \cfrac{1}{b_1 - \cfrac{1}{b_2 - \genfrac{}{}{0pt}{0}{}{\ddots}}} \end{equation} 
where $b_i \geq 2$ for $i>0$. The sequence $\{b_i\}_{i\in\mathbb{N}}$ is periodic from $b_1$, and let $b_1,\ldots,b_r$ be the minimal period of the sequence. We first define 

\begin{equation} m(p) := \frac{1}{3} \sum_{i=1}^r (b_i - 3)\,. \end{equation}

Our main theorem is the following.

\begin{theorem} \label{mainthm} If $p\equiv 3\pmod 4$ is a prime, then
\begin{equation*} h(-p) \equiv h(p) m(p) \pmod{16}\,.\end{equation*}
\end{theorem}

\begin{remark}
Combining this with Theorem \ref{Willthm}, we obtain the following congruences for $p>3$.
\begin{eqnarray*}
m(p) &\equiv& U + 2 \pmod 4\,, \\
m(p) &\equiv& 2 + pU - 2\left(\frac{T}{U}\right) \pmod 8\,.
\end{eqnarray*}
\end{remark}

This theorem is equivalent to a conjecture of Richard Guy.

\begin{corollary}[Guy's Conjecture, \cite{Guy}, \cite{Weir}]
If $p \equiv 3 \pmod 4$ is a prime, then
\begin{equation*} 
h(p) h(-p) - m(p) \equiv 
\begin{cases}
0 \, \pmod{16}\,, & \text{if } h(p) \equiv 1 \text{ or } 7\, \pmod{8} \\
8 \, \pmod{16}\,, & \text{if } h(p) \equiv 3 \text{ or } 5\, \pmod{8} 
\end{cases} \,.
\end{equation*}
\end{corollary}

To prove Theorem \ref{mainthm} we use results by Zagier \cite{Z} relating class numbers to continued fractions. Our proof also relies on properties of genus characters and Dedekind sums. It is interesting to note that the actual structure of the class group is not vital for our work.

This paper is organized as follows. In Section \ref{sectprelim}, we define our notation, give an overview of Zagier's results, and review properties of Dedekind sums and Jacobi symbols. In Section \ref{sectthm}, we prove Theorem \ref{mainthm}. Finally, in Section \ref{sectexamples}, we illustrate Theorem \ref{mainthm} using some examples.

\section{Preliminaries} \label{sectprelim}

\subsection{Some algebraic number theory}

We shall assume throughout that $p \equiv 3 \pmod 4$ is prime.  Let $\mathcal{O}_p$ denote the ring of integers of the field $\mathbb{Q}(\sqrt{p})$, and let $N(I)$ to be the norm of an ideal $I$.  It is well known (e.g., see \cite{Neukirch}) that the discriminant of $\mathbb{Q}(\sqrt{p})$ is $4p$ and that, from the theory of Pell equations, if $d,c>0$ be the smallest integer solution of Pell's equation:
\begin{equation} d^2 - p c^2 = 1\, , \label{eqnpell}\end{equation}
then $\epsilon = d + c\sqrt{p}$ is a fundamental unit of $\mathbb{Q}(\sqrt{p})$.

To prove Theorem \ref{mainthm} we make use of early work of Zagier.  Here we briefly recall his results \cite{Z}.  In \cite{Z}, $h(-p)$ is expressed using genus characters, which we describe below. 

\begin{definition}[\cite{S}] \label{defgenus}
Let $d$ be the discriminant of an algebraic number field $K$ with the decomposition $d = d_1d_2$, where each $d_i$ is also the discriminant of an algebraic number field. The \emph{genus character} $\chi$ is defined on prime ideals $\mathfrak{p}\subset \mathcal{O}_p$ not dividing $(d)$ by
\begin{equation} \chi(\mathfrak{p}) = \left(\frac{d_1}{N(\mathfrak{p})}\right) = \left(\frac{d_2}{N(\mathfrak{p})}\right)\,. \label{chip}\end{equation} 
If $\mathfrak{p} | (d)$, then one of $\left(\frac{d_1}{N(\mathfrak{p})}\right)$ and $\left(\frac{d_1}{N(\mathfrak{p})}\right)$ is zero and the other is nonzero. We take $\chi(\mathfrak{p})$ to be the nonzero value. For an arbitrary ideal $J = \prod_i\mathfrak{p}_i^{k_i}$ of $K$, we define $\chi(J) = \prod_i(\chi(\mathfrak{p}_i))^{k_i}$.
\end{definition}

\begin{remark}
It follows from the definition that $\chi(ab) = \chi(a)\chi(b)$ for any two ideals $a, b$ of $K$. For primes $\mathfrak{p}\nmid (d)$, the two definitions of the character in (\ref{chip}) will always agree.
\end{remark}

In the case of the field $\mathbb{Q}(\sqrt{p})$ with $p \equiv 3$ (mod 4), we have $d = 4p = (-4)(-p)$, and we use this factorization to determine our character. In particular, if $N(I)$ is odd for some ideal $I$, then 
\begin{equation} \chi(I)=\left(\frac{-4}{N(I)}\right)= \begin{cases} 1 & \text{if }N(I)\equiv 1 \pmod 4 \\ -1 & \text{if }N(I)\equiv 3 \pmod 4 \end{cases}\,. \end{equation}

Dedekind sums will be useful in our proofs.

\begin{definition}[\cite{RG}] 
Let $h,k$ be relatively prime integers such that $k\geq 1$. Define the \emph{Dedekind sum} $s(h,k)$ as:
\begin{equation} s(h, k) := \sum_{n = 1}^{k} \left(\left( \frac{hn}{k}\right)\right) \left(\left(\frac{n}{k}\right)\right)\, , \label{defdsum}\end{equation}
where the symbol $((x))$ is given by
\begin{equation}
((x)) := \begin{cases}
x-\lfloor x\rfloor - \frac{1}{2} & \text{if } x\not\in\mathbb{Z} \\
0 & \text{if }x\in\mathbb{Z} \end{cases}\,, \end{equation}
where $\lfloor x\rfloor$ denotes the greatest integer not exceeding $x$.
\end{definition}

It follows from the definition of the Dedekind sum that we have the following fact.
\begin{prop}[\cite{RG}]\label{ded}
Let $h'$ be an integer such that $hh'\equiv 1 \pmod{k}$. Then
\begin{equation*} s(h',k) = s(h,k)\,. \end{equation*}
\end{prop}

We will use the following key relationship between the Dedekind sum and the Jacobi symbol:

\begin{prop}[\cite{Z}]\label{dj}
If $c$ is odd, then 
\begin{equation*}\left(\frac{d}{c}\right) = (-1)^{\frac{1}{2}\left(\frac{c - 1}{2} - 6c s(d,c)\right)}\,. \end{equation*}
\end{prop}

\subsection{Overview of Zagier's work}
In \cite{Z}, Zagier relates class numbers of $\mathbb{Q}(\sqrt{-p})$ with the negative continued fraction expansion of $\sqrt{p}$.  Here, we describe some of the main ideas relating to this result which we use in our work.  A particular subset of $SL(2,\Z)$ known as the hyperbolic matrices is important in this connection.

\begin{definition}
Let $A = \begin{pmatrix} a & b \\ c & d \end{pmatrix}  \in SL(2, \mathbb{Z})$. $A$ is a \emph{hyperbolic} matrix if $|\textup{tr}(A)|>2$.  
\end{definition}

A key idea in Zagier's work is the one-to-one correspondence \begin{equation} A \longleftrightarrow (M, V) \end{equation} between conjugacy classes of hyperbolic matrices $A$ and equivalence classes of pairs $(M, V)$, where $M$ is a rank $2$ free $\mathbb{Z}$-module inside a real quadratic field $K$, and $V$ is a free abelian group of rank $1$ generated by a totally positive unit $\epsilon \in K$ satisfying $\epsilon M=M$. Two pairs $(M,V)$ and $(M',V')$ are equivalent if $M'=\alpha M$ for some $\alpha\neq 0$ in $K$ and $V'=V$.

Given such a pair $(M,V)$, let $\{\beta_1,\beta_2\}$ be an {\em oriented basis} of $M$ in the sense that $\beta_2>0$ and $\beta_1\beta_2' - \beta_1'\beta_2 >0$, where $\beta_i'$ denotes the conjugate radical of $\beta_i$ for $i=1,2$. Let $\epsilon$ be a generator of $V$. We get a corresponding matrix $A= \begin{pmatrix} a & b \\ c & d \end{pmatrix} \in SL(2,\mathbb{Z})$ by
\begin{eqnarray*}
\epsilon\beta_1 &=& a\beta_1 + b\beta_2\,, \\
\epsilon\beta_2 &=& c\beta_1 + d\beta_2\,. \end{eqnarray*}

Conversely, given a hyperbolic matrix $A\in SL(2,\mathbb{Z})$, let $w$ be the larger of its two fixed points under its action on $\C$ by linear fractional transformation.  We set $\beta_1=w$, $\beta_2=1$, $M$ to be the module generated by $w,1$ and $\epsilon>1>\epsilon'>0$ the eigenvalues of $A$.\\

For our purposes, we will generally let $M$ be some ideal of the ring of integers $\mathcal{O}_p$ and $V$ be the (free rank $1$) group of totally positive units of $\mathcal{O}_p$.

Given any hyperbolic matrix $A$ with $c>0$, it is known from \cite{Z} that $A$ is conjugate in $SL(2,\mathbb{Z})$ to a product of the form

\begin{equation}  A = \begin{pmatrix} b_1 & -1 \\ 1 & 0 \end{pmatrix} \cdots \begin{pmatrix} b_r & -1 \\ 1 & 0 \end{pmatrix} \label{matrixbform},\end{equation} 

\noindent with $b_1, \ldots, b_r \in \mathbb{Z}$ with $b_i\geq 2$ and at least one $b_i>2$.  The sequence $b_1,\ldots, b_r$ is unique up to cyclic permutation, and can be computed using the fixed points of $A$ viewed as acting on $\mathbb{C}$ via linear fractional transformation.  As before, let $w$ be the larger solution to the equation $Az = z$. Since $w$ satisfies a quadratic equation over $\mathbb{Q}$, its negative continued fraction expansion eventually becomes periodic. If $V$ is the group of all totally positive units which leave $M$ invariant for the pair $(M,V)$ corresponding to $A$, then $b_1, \ldots b_r$ is the period of this sequence.  This is clearly true if $M$ is an ideal of $\mathcal{O}_p$ and $V$ is taken to be the full group of totally positive units of $\mathcal{O}_p$.\\

Given $A = \begin{pmatrix} a & b \\ c & d \end{pmatrix} \in SL(2, \mathbb{Z})$, with $c>0$, we define the integer
\begin{equation} n_A = \frac{a + d}{c} - 3 - 12s(d,c)\,,\label{na}\end{equation} 
where $s(d,c)$ is the Dedekind sum defined in (\ref{defdsum}).

If we write $A$ in the form in (\ref{matrixbform}), then Zagier shows that (\ref{na}) simplifies to
\begin{equation} n_A = \sum_{i=1}^r \left(b_i-3\right)\, , \end{equation} in \cite[p.90]{Z}.

\begin{remark} For the purpose of this paper, it suffices to treat (\ref{na}) as a definition. Zagier \cite{Z} defines $n_A$ in terms of a modular form and proves that $n_A$ can be computed as in (\ref{na}).
\end{remark}

We now state one of Zagier's results.

\begin{theorem}[\cite{Z}] \label{thmz2}
If $p \equiv 3 \pmod{4}$, then we have
\begin{equation*} h(-p) = \frac{1}{3}\sum_{C} \chi(I)n(I), \end{equation*}
where the sum runs through the ideal classes $C$ of $\mathbb{Q}(\sqrt{p})$, $\chi$ is a genus character as defined in Definition \ref{defgenus} and $n(I) = n_A$, where $A$ is the matrix corresponding to the action of the fundamental unit of $\mathbb{Q}(\sqrt{p})$ on a basis of $I$.
\end{theorem}

\begin{remark}
As conjugation by an element of $SL(2,\mathbb{Z})$ does not change $n_A$, it does not matter which basis of $I$ is chosen.
\end{remark}

In the case of the quadratic number field $\mathbb{Q}(\sqrt{p})$, the unit ideal of $\mathcal{O}_p$ with basis $\{\sqrt{p},1\}$ corresponds to the matrix $A= \begin{pmatrix} d & cp \\ c & d \end{pmatrix}$, where $d + c\sqrt{p}$ is the fundamental totally positive unit. The largest fixed point of $A$ is $\sqrt{p}$, which we can write as a negative continued fraction as in (\ref{contfracp}), and we have the formula for $n_A$ as in (\ref{matrixbform}). This implies the following theorem.
\begin{theorem}[\cite{Z}] \label{thmz}
Let $p \equiv 3 \pmod{4}$, and assume that $h(p) = 1$. Let $b_1, \ldots, b_r$ be the period of the negative continued fraction expansion of $\sqrt{p}$ as in (\ref{contfracp}). Then we have 
\begin{equation*} h(-p) = \frac{1}{3}\sum_{i = 1}^{r} (b_i - 3)\,. \end{equation*}
\end{theorem}

\section{Proof of Theorem \ref{mainthm}} \label{sectthm}
Let $p\equiv 3\pmod 4$ be a prime and let $\mathcal{O}_p$ be the ring of integers of $\mathbb{Q}(\sqrt{p})$. Let $\chi(I)$ and $n(I)$ be as defined in Section \ref{sectprelim}.  

We have by Theorem \ref{thmz2} that 
\begin{equation}\label{h(-p)} h(-p) = \frac{1}{3}\sum_C \chi(I)n(I)\,.
\end{equation}

Since $\chi(I)n(I)$ is independent of the representative $I\in C$ chosen, we may define 
\begin{equation} t_C := \chi(I)n(I)\,. \end{equation}  We shall first make an observation about choosing an appropriate representative of each ideal class $C$.  This will allow us to prove two lemmas about the $t_C$, from which the theorem then follows.

\begin{observation}For any ideal class $C$, there exists a representative ideal $I\in C$ such that $\chi(I)=1$, $N(I)$ is odd, and $I = (a+\sqrt{p},b)$ for integers $a$ and $b$.  
\end{observation}
\begin{proof}
Take any ideal $I\in C$.  We first show that we can choose $I$ such that $N(I)$ is odd.

If $4|N(I)$, then for any $a+b\sqrt{p}\in I$, we have $4|a^2-pb^2$ and since $p\equiv 3\pmod{4}$, we have that $2|a,b$ and so $I = (2)I'$ for some ideal $I'\in C$.  

If $N(I) \equiv 2\pmod{4}$, then $N((1+\sqrt{p})I)\equiv 4\pmod{8}$ and so $(1+\sqrt{p})I = (2) I'$ where $I'\in C$ has odd norm.  Thus, we may choose $I$ to have odd norm.  

If $\chi(I) = -1$, then note that $(\sqrt{p})$ has odd norm and $\chi((\sqrt{p}))=-1$, so $I' = (\sqrt{p})I$ has character $1$ and odd norm.

We also note that any ideal class $C\in \mathcal{O}_p$ contains a representative of the form $(x+\sqrt{p}, y)$ for integers $x,y$.  Indeed, given a representative ideal in the form $I=(a+b\sqrt{p},c+d\sqrt{p})$, we can apply the Euclidean algorithm to obtain the generators $I= (a'+b'\sqrt{p}, c')$ where $c'$ is the smallest integer in $I$.  If $(b',c')=1$, then we can immediately reduce $I$ to the desired form.  Otherwise, suppose that $q|b'$ and $q|c'$ for some prime $q$.  Then, since $a'^2 - pb'^2$ is an integer in $I$, we have by definition that $c'|a'^2-pb'^2$.  It follows that $q|a'^2-pb'^2$ and hence $q|a'$.  Therefore, we may write $I = (q)I'$ for some ideal $I'\in C$.  Repeating the process, we find the desired ideal.  

Note that since we only factor out ideals of the form $(q)$, this process will not change $\chi(I)$ or the parity of $N(I)$.  Therefore, we may, for any ideal class, choose a representative $I = (a+\sqrt{p},b)$ such that $\chi(I) = 1$ and $N(I)$ is odd.  
\end{proof}

\begin{lemma}\label{inverses}
We have that $t_C = t_{C^{-1}}$.  
\end{lemma}

\begin{proof}
Throughout, let $(x,y)=(d,c)$ be the minimal integer solution to the Pell equation $x^2-py^2 = 1$.  Let $I= (a+\sqrt{p},b)\in C$ be a representative chosen as above, and let $J = (-a+\sqrt{p}, b)\in C^{-1}$ be a representative, where the generators have been oriented.  Note that $\chi(I)=\chi(J)$ because these ideals have the same norm.  Thus, it suffices to show $n(I)=n(J)$.  

The relevant matrix for computation of $n(I)$ is given by a change of basis from $(\sqrt{p},1)$ to $(a+\sqrt{p},b)$ as follows:

\[ \left( \begin{array}{cc}
1 & a  \\
0 & b \end{array} \right)  
\left( \begin{array}{cc}
d & cp  \\
c & d \end{array} \right)  
\left( \begin{array}{cc}
1 & -\frac{a}{b}  \\
0 & \frac{1}{b} \end{array} \right)  = 
\left( \begin{array}{cc}
d+ac & \frac{c(p-a^2)}{b}  \\
cb & d-ac \end{array} \right).
\]

Similarly, we have the following matrix which computes $n(J)$:

\[ \left( \begin{array}{cc}
1 & -a  \\
0 & b \end{array} \right)  
\left( \begin{array}{cc}
d & cp  \\
c & d \end{array} \right)  
\left( \begin{array}{cc}
1 & \frac{a}{b}  \\
0 & \frac{1}{b} \end{array} \right)  = 
\left( \begin{array}{cc}
d-ac & \frac{c(p-a^2)}{b}  \\
cb & d+ac \end{array} \right).
\]

Finally, by (\ref{na}), we simply need to show that $s(d+ac,cb) = s(d-ac,cb)$.  Observe that 

\begin{equation*}
(d+ac)(d-ac) = d^2 - a^2c^2 = d^2-pc^2 + (p-a^2)c^2 \equiv 1\pmod{cb},
\end{equation*}
because $d^2-pc^2 =1$ and $b$, being the norm of $I$, divides the norm of $a+\sqrt{p}$.  Hence, $d+ac$ and $d-ac$ are inverses modulo $cb$ and the conclusion follows directly from Proposition \ref{ded}.  
\end{proof}

\begin{lemma}\label{equalmod8}
For any two classes $C,C'$, we have that $t_C\equiv t_{C'}\pmod{8}$.   
\end{lemma}
\begin{proof}
It suffices to show the statement when $C'$ is the identity.  Let $I = (a+\sqrt{p},b)$ be such that $N(I)= b$ is odd and $\chi(I) = 1$.  We compute $n((1))$ with respect to the oriented basis $(a+\sqrt{p},1)$ of the unit ideal.  Then, as before, we may compute $n((1))$ with the matrix:
\[ \left( \begin{array}{cc}
1 & a  \\
0 & 1 \end{array} \right)  
\left( \begin{array}{cc}
d & cp  \\
c & d \end{array} \right)  
\left( \begin{array}{cc}
1 & -a  \\
0  & 1 \end{array} \right)  = 
\left( \begin{array}{cc}
d+ac & c(p-a^2) \\
c & d-ac \end{array} \right),
\]

\noindent and we may compute $n(I)$ with the matrix
\[ \left( \begin{array}{cc}
d+ac & \frac{c(p-a^2)}{b}  \\
cb & d-ac \end{array} \right).
\]
as found in the previous lemma.  

Therefore, by (\ref{na}) it suffices to show that 
\begin{equation*}
\frac{2d}{c} - 12s(d-ac, c) \equiv \frac{2d}{cb} - 12s(d-ac,cb) \pmod{8}\,.  \end{equation*}
However, $b$ was chosen to be odd and $d$ is even by \cite[Theorem 1.1]{ZY14}.  Thus, it suffices to show that $12s(d-ac,c)\equiv 12s(d-ac,cb)\pmod{8}$.  However, note that by (\ref{dj}) we have the following equations:
\begin{eqnarray}
\left( \frac{d-ac}{cb}\right) &=& (-1)^{\frac{1}{2}(\frac{cb-1}{2} -6cbs(d-ac,cb))} \\
\left( \frac{d-ac}{c}\right) &=& (-1)^{\frac{1}{2}(\frac{c-1}{2} -6cs(d-ac,c))}.
\end{eqnarray}

Dividing these equations, we obtain:
\begin{equation}
\left( \frac{d-ac}{b}\right) = (-1)^{\frac{c(b-1)}{4} + 3c(s(d-ac,c)-bs(d-ac,cb))}.
\end{equation}

We now note that $b\equiv 1\pmod{4}$ because $b$ is odd and is the norm of an ideal with character 1.  Furthermore, $c$ is odd because $d$ is even.  Hence, 
\begin{equation*}
(-1)^{\frac{c(b-1)}{4}} = (-1)^{\frac{b-1}{4}} = (-1)^{\frac{b^2-1}{8}} = \left(\frac{2}{b}\right).  
\end{equation*}

Finally, because $p\equiv 3\pmod{4}$, we have that $2(d-c\sqrt{p}) = (R-S\sqrt{p})^2$ for integers $R,S$ by \cite{Will82}.  Hence, $\left(\frac{2(d-ac)}{b}\right) = \left(\frac{2(d-c\sqrt{p})}{b}\right) =1$.  Therefore, $3c(s(d-ac,c)-bs(d-ac,cb))\equiv 0\pmod{2}$ and so $12c(s(d-ac,c) - s(d-ac,cb))\equiv 0\pmod{8}$ as desired.  

\end{proof}

We now deduce Theorem \ref{mainthm} from the lemmas.
\begin{proof}[Deduction of Theorem \ref{mainthm} from Lemmas \ref{inverses} and \ref{equalmod8}]

Observe that $m(p) = \frac{1}{3}t_{\text{id}}$.  Therefore, since the sum (\ref{h(-p)}) goes over $h(p)$ terms, we have that 
\begin{equation}
h(-p) - m(p)h(p) = \frac{1}{3}\left(\sum_C t_C\right) - m(p)h(p) = \sum_C \left(\frac{t_C}{3} - m(p)\right).
\end{equation}

The term corresponding to the identity is just 0.  By Lemma \ref{equalmod8}, we have that all the remaining terms are divisible by 8.  Finally, since $h(p)$ is odd (see, for example, \cite{ZY14}), Lemma \ref{inverses} implies that these remaining terms come in equal pairs modulo 16.  Hence Theorem \ref{mainthm} follows as desired.  

\end{proof}

\section{Examples} \label{sectexamples}
In this section, we illustrate Theorem \ref{mainthm} using three examples. The case where $h(p)=1$ has already been studied by Zagier \cite{Z}, and $h(-p)=m(p)$ by Theorem \ref{thmz}. We give examples of cases where $h(p)>1$. In what follows, we use the standard continued fraction notation, whereby we denote the continued fraction in (\ref{contfracp}) as $[b_0;\overline{b_1,\ldots,b_r}]$.

\subsection{Example 1: $p=79$}
The smallest prime $p\equiv 3$ (mod 4) such that $h(p)>1$ is $p=79$. We have $h(79)=3$ and $\sqrt{79}=[9;\overline{9,18}]$. Therefore, we have
\begin{equation*} m(79)=\frac{1}{3}((9-3)+(18-3))=7\,. \end{equation*}

Let $C_0$, $C_1$, and $C_2$ be the three ideal classes, with $C_0$ the class of principal ideals. Then $t_{C_0}=3m(79)=21$, and by Lemma \ref{inverses}, $t_{C_1}=t_{C_2}$. Therefore, we only need to calculate $t_{C_1}$.

We choose $I=(1+\sqrt{79},3)$ as a representative of $C_1$ (we may switch $C_1$ and $C_2$ so that this holds, if necessary). Since $N(I)=3$, we have $\chi(I)=-1$. To compute $n(I)$, we find the matrix $A_I \in SL(2,\mathbb{Z})$ corresponding to this ideal and basis,
\[A_I=\left(\begin{array}{cc} 1 & 1 \\ 0 & 3 \\\end{array} \right)\left(\begin{array}{cc} d & 79c \\ c & d \\\end{array} \right)\left(\begin{array}{cc} 1 & 1 \\ 0 & 3 \\\end{array} \right)^{-1},\]
where the fundamental unit is $\epsilon=d+c\sqrt{79}$. This comes from taking the action of $\epsilon$ with respect to the basis $(\sqrt{79},1)$ and conjugating it to change the basis to that of $I$.

To compute $n(I)=n_{A_I}$, we view $A_I$ as the corresponding fractional linear transformation, take its largest fixed point, and compute the negative continued fraction expansion for that quadratic irrational. Now the advantage of writing $A_I$ in the above form is clear: since the largest fixed point of $\left(\begin{array}{cc} d & 79c \\ c & d \\\end{array} \right)$ is $\sqrt{79}$, the largest fixed point of $A_I$ is
\[\left(\begin{array}{cc} 1 & 1 \\ 0 & 3 \\\end{array} \right)(\sqrt{79})=\frac{1+\sqrt{79}}{3}\,.\]
The negative continued fraction expansion of $\frac{1+\sqrt{79}}{3}$ is $[4;\overline{2,2,4,3,7}]$, so
\[n(I)=(2-3)+(2-3)+(4-3)+(3-3)+(7-3)=3\,.\]
Therefore, $t_{C_1}=t_{C_2}=-3$, so the form of our expression is 
\begin{equation*} h(-79)=\frac{1}{3}(21-3-3)=5\,. \end{equation*}

\subsection{Example 2: $p=439$}
The smallest $p$ such that $h(p)>3$ is $p=439$, and $h(439)=5$. Since $\sqrt{439}=[21;\overline{21,42}]$, $m(439)=19$. Similarly to the previous example, let $C_0$, $C_1$, $C_2$, $C_3$, and $C_4$  be the five ideal classes, where we label them such that the cyclic group structure on the ideal classes is the same as that on their indices. Then $t_{C_0}=3m(439)=57$, $t_{C_1}=t_{C_4}$, and $t_{C_2}=t_{C_3}$.

Now choose $I=(7+\sqrt{439},13)$ to be a representative of $C_1$. Similarly to the previous example, we compute the negative continued fraction expansion of $\frac{7+\sqrt{439}}{13}$. This is $[3;\overline{2,2,2,2,2,2,2,2,2,2,2,2,3,3,2,2,2,2,2,5}]$, so $n(I)=-15$. Now, $N(I)=13\equiv 1$ (mod 4), so $\chi(I)=1$. Therefore, $t_{C_1}=t_{C_4}=-15$.

We next take $I=(13+\sqrt{439},18)$, which is a representative of $C_2$. Correspondingly, $\frac{13+\sqrt{439}}{18}$ has a negative continued fraction expansion $[2;\overline{9,5,5,2,3}]$. Therefore, $n(I)=9$. Finally, $N(I)=18$, and since $439\equiv 7$ (mod 8), $\chi(I)=\left(\frac{-439}{18}\right)=1$. Therefore, $t_{C_2}=t_{C_3}=9$, so here our expression will be of the form $h(-439)=\frac{1}{3}(57-15+9+9-15)=15$.

\subsection{Example 3: $p=43063$}
For our final example, we take $p=43063$, which is the smallest prime congruent to $3$ modulo $4$ such that the corresponding class group is not cyclic; it is isomorphic to $\mathbb{Z}/3\mathbb{Z}\times\mathbb{Z}/3\mathbb{Z}$. Computation via Sage yields $m(43063)=193$ and the following table of terms $t_C$, arranged by the powers of two generators of the class group that represent that class:
\begin{center}
\begin{tabular}{c | c | c | c |}
  & 0 & 1 & 2 \\ \hline
0 & 579 & -21 & -21 \\ \hline
1 & 51 & -141 & -69 \\ \hline
2 & 51 & -69 & -141 \\ \hline
\end{tabular}
\end{center}
yielding $h(-43063)=73$.

\section{Acknowledgments}
We would like to thank Ken Ono for his advice and encouragement. This paper was written while all the authors were participants in the 2014 Emory Math REU, and as such we would like to thank the NSF for its support and the Emory Department of Mathematics and Computer Science for its hospitality.

\end{document}